\documentclass[11pt]{article}
\vspace{7cm} \textwidth 14cm \textheight 21.6cm
\usepackage{amsmath}
\usepackage{amsfonts}
\usepackage{latexsym,amsthm,amscd}
\newtheorem{thm}{Theorem}

\newtheorem{prop}{Proposition}

\newtheorem{lem}{Lemma}
\newtheorem{Def}{Definition}

\newcounter{alphthm}
\newcommand{\be}{\begin{equation}}
\newcommand{\ee}{\end{equation}}
\newcommand{\ben}{\begin{enumerate}}
\newcommand{\een}{\end{enumerate}}
\newcommand{\beq}{\begin{eqnarray}}
\newcommand{\eeq}{\end{eqnarray}}
\newcommand{\beqn}{\begin{eqnarray*}}
\newcommand{\eeqn}{\end{eqnarray*}}

\newcommand{\bpf}{\begin{proof}}
\newcommand{\epf}{\end{proof}}
\newcommand{\bl}{\begin{lem}}
\newcommand{\el}{\end{lem}}
\newcommand{\bp}{\begin{prop}}
\newcommand{\ep}{\end{prop}}
\newcommand{\bd}{\begin{Def}}
\newcommand{\ed}{\end{Def}}
\newcommand{\bt}{\begin{thm}}
\newcommand{\et}{\end{thm}}

\title{Classification of complete Finsler manifolds through a second order differential equation}
\author{\small  A. ASANJARANI and B. BIDABAD\\ }
\date{\small{ Faculty of
Mathematics and Computer Sciences, Tehran Polytechnic (Amirkabir
University of Technology), Hafez Ave., 15914 Tehran, Iran.}}

\begin{document}
\maketitle
--------------------------------------------------------------------------------------------------
\section*{Abstract}
\hspace{.5 cm}By using a certain second order differential equation,
the notion of adapted coordinates on Finsler manifolds is defined
and some classifications of complete Finsler manifolds are found.
Some examples of  Finsler metrics, with positive constant sectional
curvature, not
 necessarily of Randers type nor projectively flat, are found.
 This work generalizes
some results in Riemannian geometry and open up, a vast area of
research on Finsler geometry.\\\\
\emph{ MSC: Primary 53C60; Secondary 58B20.}\\\\
 \emph{Keywords }: Finsler, conformal, constant curvature, second order
differential equation.\\
---------------------------------------------------------------------------------------------------
\section*{Introduction}
\setcounter{equation}{0}
 \hspace{.5 cm}Differential equations play
an essential role in the study of global differential geometry,
particularly the second order differential equation
 $\nabla\nabla\rho= \phi g$, which
 appears often in the study of a Riemannian manifold $(M,g)$. In the above
equation,  $\rho$ is the divergence of a  conformal Killing vector
field,  $\phi$ is a $C^\infty$ function on $M$ and $\nabla$ is the
covariant derivative of Levi-Civita connection. This differential
equation has proved to be  very fruitful and has been studied by
many authors  not to be mentioned here, see for instance [10], [17],
[20], [21] and [22]. Geometrically, the existence of a solution to
this differential equation is equivalent to the existence of a
certain conformal transformation,
 which takes  geodesic circles into geodesic circles -
a geodesic circle being a curve with constant curvature and zero
torsion  and is a generalization of  circle in the Euclidean space -
see for example [11] and [21].  In Physics, this differential
equation is closely connected to the study of collineations in
General Relativity [15]. The above differential equation appears
also in the study of pseudo-Riemannian
manifolds, see for example [8] and [14].\\
In Finsler geometry, this differential equation has been
investigated in [1] and [2], using a method of calculus of
variation. The results  there obtained may be considered as a very
special case of the main theorem  presented in this work.

We propose to consider some possible applications of this
differential equation in  a Finslerian setting, having in mind the
following remarks: First of all  existence of a solution permits the
definition of  a new adapted coordinate system, which will somehow
play the same role in Finsler geometry try as  the normal coordinate
system in Riemannian geometry($\S$1, $\S$2 and $\S$3) in which it
has proved to be a powerful tool, while its usefulness has been
limited in the Finsler geometry. In fact, in the latter case the
exponential map is only $C^1$ at the zero section of $TM$, while it
is $C^\infty$ in the former  case [1], [3]. Next,
 the following classification theorem can be proved($\S$5).\\
\textit{\textbf{Theorem:}}\emph{Let $(M,g)$ be a connected complete
Finsler manifold of dimension $n\geq 2$. If $M$ admits a non-trivial
solution of $ \nabla^{H}\nabla^{H}\rho= \phi g$, where $\nabla^{H}$
is the Cartan horizontal covariant derivative,
 then  depending on the number of critical points of $\rho$, i.e. zero, one or two respectively,  it is conformal to\\
 \textbf{(a)}   A direct product $J\times\overline{M}$ of an open
interval $J$ of the real line and an       $(n-1)$-dimensional
complete
Finsler manifold $\overline{M}$.\\ \textbf{(b)} An n-dimensional Euclidean space.\\
 \textbf{(c)} An n-dimensional  unit sphere in an Euclidean
 space.}\\
It should be remarked  that the role played by  Cartan derivative in
the above theorem is essential and can not be replaced by, for
example the Berwald derivative. Well known examples of Finsler
metrics with positive constant curvature are either of Randers type
or are projectively flat [3],[4],[5],[7],[9]. As yet another next
application of this differential equation, we find some examples of
 Finsler metrics with positive constant sectional curvature which are not
 necessarily of Randers type nor projectively flat($\S$4).
More precisely, by using a Finsler metric with positive constant
sectional curvature together with the above theorem, one can
construct Finsler metrics of positive constant
 curvature in higher dimensions($\S$6).
\section*{Preliminaries.}
  Let $M$ be a real n-dimensional  manifold of class $C^ \infty$. We
denote by  $TM\rightarrow M$ the  bundle
  of tangent vectors and by $ \pi:TM_{0}\rightarrow M$ the fiber bundle
of non-zero tangent vectors.
  A {\it{Finsler structure}} on $M$ is a function
$F:TM \rightarrow [0,\infty )$, with the following properties: (I)
$F$ is differentiable ($C^ \infty$) on $TM_{0}$; (II) $F$ is
positively homogeneous of degree one in $y$, i.e.
 $F(x,\lambda y)=\lambda F(x,y),  \forall\lambda>0$, where $(x,y)$
 is an element of $TM$.
(III) The Hessian matrix of $F^{2}$, $(g_{ij}):=\left({1 \over 2}
\left[ \frac{\partial^{2}}{\partial y^{i}\partial y^{j}} F^2
\right]\right)$, is positive definite on $TM_{0}$.
 A \textit{Finsler
manifold} is a pair consisting of a differentiable manifold $M$ and
a Finsler structure $F$ on $M$. The tensor field $g$ with the
components $g_{ij}$ is called the  \emph{ Finsler metric tensor}.
Hereafter, we denote a Finsler manifold by $(M,g)$. Let
$V_vTM=ker\pi_{*}$ be the set of vectors tangent to the fiber
through $v\in TM_0$. Then a \emph{vertical vector bundle} on $M$ is
defined by $VTM := \bigcup_{_{v\in TM_0}}V_vTM$. A
\textit{non-linear connection} on $TM_0$ is a complementary
distribution $HTM$  for $VTM$ on $TTM_0$. Therefore we have the
decomposition
 $
TTM_0 =VTM\oplus HTM
 $.
 The pair $(HTM,\nabla)$, where  $HTM$ is a
non- linear connection on $TM$ and $\nabla$ a linear
  connection on $VTM$,  is called a
  \textit{Finsler connection} on the manifold $M$.
Using the local coordinates $(x^{i},y^{i})$ on $TM$, called the line
elements, we have the local field of frames
$\{\frac{\partial}{\partial x_{i}},\frac{\partial}{\partial
y_{i}}\}$ on $TTM$.
  Given a non-linear connection one can choose a local field of frames
$\{\frac{\delta}{\delta x_{i}},\frac{\partial}
  {\partial y_{i}}\}$ adapted to the above decomposition i.e.
$\frac{\delta}{\delta x_{i}}\in {\Gamma}(HTM)$
   and $\frac{\partial}{\partial y_{i}}\in {\Gamma}(VTM)$, the set of
vector fields on $HTM$
  and $VTM$ respectively. Here
  $
  \frac{\delta}{\delta x_{i}}:=\frac{\partial}{\partial
x_{i}}-G^{j}_{i}\frac{\partial}
  {\partial y_{j}},
 $
and  $G^{j}_{i}(x^{},y^{})$  are coefficients of the non-linear
connection.
 A \textit{Cartan connection} is a metric
Finsler connection and its  coefficients are defined by
$\Gamma^{\ast{}i}_{jk}=\frac{1}{2}g^{ih}(\frac{\delta g_{hk}}{\delta
x^{j}}+ \frac{\delta g_{hj}}{\delta x^{k}}- \frac{\delta
g_{jk}}{\delta x^{h}}) $.
 Also we can write $\Gamma^{\ast{}i}_{jk}= \gamma^{i}_{jk}-C^{i}_{jr}G^{r}_{k}+
C^{r}_{jk}G^{i}_{r}- C^{i}_{kr}G^{r}_{j},$
 where  $\gamma^{i}_{jk}$ are \emph{formal Christoffel symbols} of the
second kind given by
 $
\gamma^{i}_{jk}=\frac{1}{2}g^{ih}(\frac{\partial g_{hk}}{\partial
x^{j}}+ \frac{\partial g_{hj}}{\partial x^{k}}- \frac{\partial
g_{jk}}{\partial x^{h}})$.
   Let
$\overline{M}$ and  $M$ be two differentiable manifolds of dimension
$m$ and $m+n$  respectively and let   $(u^{\alpha})$ and $(x^{i})$
be
  local coordinate systems  on them.
  We denote by  $(u^\alpha ,
v^\alpha)$ and  $(x^i, y^i)$ pairs of position and direction of the
line elements of $T\overline{M}$ and $TM$, where $\alpha, \beta,
...$ and $i, j, ...$ run over the range $1,...,m$ and $1,..., m+n$
respectively. Let $f:\overline{M}\rightarrow M$ be a smooth map,
given by
 $
(u^{1}, ..., u^{m})\rightarrow x^{i}(u^{1}, ..., u^{m})$. The
differential mapping of $f$ is
$$ \begin{array}{lcl}
\qquad{f_{*}:T_{u}\overline{M}\rightarrow T_{x}M },\\
(u^{\alpha}, v^{\alpha})\rightarrow (x^{i}(u),  y^{i}(u,v)),
 \end{array} $$
where $y^{i}(u,v)= B_{\alpha}^{i}v^{\alpha}$ and
$B_{\alpha}^{i}=\frac{\partial x^{i}}{\partial u^{\alpha}}$. If
$f_{*}$ is injective at every point $u$ of $\overline{M}$, that is,
if rank $[B^{i}_{\alpha}]=m$, then $\overline{M}$ is called an
\emph{immersed submanifold} or simply  a \emph{submanifold} of $M$.
Next, consider an $(m+n)$-dimensional Finsler manifold $(M,g)$. The
Finsler structure $F$ induces on $T\overline{M}$ a Finsler structure
$\overline{F}$ defined by $ \overline{F}(u,v):= F(x(u), y(u,v)).$
Putting $\overline{g}_{\alpha\beta}:={1 \over 2}
 \frac{\partial^{2} \overline{F}^2}{\partial v^{\alpha}\partial
v^{\beta}} $, one obtains by direct calculation $
\overline{g}_{\alpha\beta}(u,v)= g_{ij}(x(u),
y(u,v))B_{\alpha\beta}^{ij},$
 where $B_{\alpha\beta}^{ij}=B_{\alpha}^{i}B_{\beta}^{j}$. Therefore
the pair $(\overline{M},\overline{g})$ is a Finsler manifold,
called \emph{Finsler submanifold } of $(M,g)$. \\
A diffeomorphism $f : (M,g)\rightarrow (N,h)$ between
$n$-dimensional
 Finsler manifolds $(M,g)$ and $(N,h)$
is called \textit{conformal}  if each $(f_{*})_p $ for $p\in M$ is
angle-preserving, and in this case two Finsler manifolds are called
\textit{conformally equivalent} or simply \textit{conformal}. If
$M=N$ then$f$ is called a \textit{conformal transformation}. It can
be easily checked that a diffeomorphism $f$ is conformal if and only
if, $f^{*}h = \rho g$ for some positive function $\rho: M\rightarrow
I\!\!R^+$. The diffeomorphism $f$ is called an \emph{homothety} if
$\rho=$constant and an \textit{isometry} if $\rho=1$.
 Now let's consider two Finsler manifolds $(M,g)$ and
$(M,\overline{g})$, then these two manifolds are conformal if and
only if, $\overline{g} = \rho(x)\ g$.

Throughout this paper, all manifolds are supposed to be connected.
\section{ Finsler manifolds admitting a non-trivial solution of $\nabla^{H}\nabla^{H}\rho= \phi g$.}
 \setcounter{equation}{0}
  Let $(M,g)$ be an n-dimensional Finsler manifold and   $\rho:M\rightarrow [0,\infty)$  a scalar  function on $M$
   given by the
following second order differential equation
 \be
 \label{c-field}
\nabla^{H}\nabla^{H}\rho= \phi g,
 \ee
where $\nabla^{H}$ is the Cartan horizontal covariant derivative and
$\phi$ is a function of $x$ alone, then we say that the Eq.
(\ref{c-field}) has a solution $\rho$.\footnote{In the local
coordinate system, the Eq. (\ref{c-field}) is given by
 $
\nabla^{H}_k\rho_l= \frac{\partial}{\partial x^k}\rho_l +
\Gamma^{\ast{}j}_{lk}\rho_j= \phi g_{kl}
 $, where $\rho _l =\frac{\partial \rho}{\partial x^l}$.}

In this section we consider the non-trivial (i.e. non-constant)
solution $\rho$
 of the
 Eq.(\ref{c-field}). The
connected component of a regular hypersurface defined by
$\rho=constant$, is called a\emph{ level set of $\rho$}. Let's
denote by $\verb"grad" \rho$ the gradient vector field of $ \rho$
which is locally
 written in the form $\verb"grad" \rho =  \rho^{i}\frac{\partial}{\partial
x^i}$, where $\rho^i = g^{ij} \rho_j$ and $i, j, ...$ run over the
range $1,...,n$. Contracting (\ref{c-field}) with $\rho^{k}$, we get
$ \rho^{k}(\nabla_{k}^{H}\rho^{l})=\phi\rho^{l}$ or equivalently $
\rho^k\frac{\partial \rho^l }{\partial
x^k}+\Gamma^{\ast{}l}_{jk}\rho^j\rho^k=\phi\rho ^{l},
 $ which shows  that the trajectories of the vector field $\verb"grad" \rho
$ are geodesic arcs.

 Therefore we can choose a local coordinates
  $(u^1=t, u^2,...,u^n)$ on $M$  such that $t$ is  parameter of the geodesic containing a trajectory of the
vector field $\verb"grad" \rho$ and the level sets of $\rho$ are
      defined by $t=$constant, called respectively the \emph{ $t$-geodesic}
      and the  \emph{t-levels} of $\rho$.
In the local coordinates $(t, u^2,...,u^n)$, $t$-levels of $\rho$
are defined by $t$=constant, so $\rho$ may be considered as a
function of $t$ alone.
 In the sequel we will refer to this coordinates as   an \emph{adapted
coordinates}.

  The differential equation  of $t$-geodesics are given by
 \be \label{2.12}
\frac{d}{dt}\frac{du^i}{dt}+ 2G^{i}=\varrho \frac{du^i}{dt},
 \ee
 where
$\varrho$ is a scalar function of $t$  and
$G^{i}=\frac{1}{2}\Gamma^{*{}i}_{jk}\frac{du^j}{dt}\frac{du^k}{dt}$
are spray coefficients.
  Since $\frac{du^{i}}{dt}=\delta_{1}^{i}$, where
$\delta_{j}^{i}$ is the
  Kronecker symbol, from the above geodesic differential equation
  we have
 \be \label{G1} G^1=\frac{1}{2}\varrho  \quad{ and}\quad{G^\alpha=0, \alpha=2,3,...,n} \ee

  \bp
Let $(M,g)$ be a Finsler manifold  admitting a non-trivial solution
of (\ref{c-field}). Then
 $g_{1\beta}=g_{\beta1}=0$, where $\beta=2,3,...,n$ and
$g_{11}=1$.
 \ep
 \bpf
 The $t$-geodesics are normal to the
$t$-levels of $\rho$, so we have
 \be \label{2.11}
 g_{1\beta}=g_{\beta1}=0, \quad{\quad{\beta=2,3,...,n}}.
 \ee
 Putting (\ref{2.11}) in the definition of
$\Gamma^{\ast}{}_{11}^{i}$ and using (\ref{G1}) we get
 \be
 \label{2.13}
\frac{1}{2}g^{i1}\delta_{1}g_{11}-\frac{1}{2}g^{i\alpha}\delta_{\alpha}g_{11}=0.
 \ee
By replacing the index $i$ by $\beta$ in the above equation, we get
$\delta_{\alpha}g_{11}=0$. As a consequence of (\ref{G1}), the
equation $
  \frac{\delta}{\delta x_{i}}=\frac{\partial}{\partial
x_{i}}-G^{j}_{i}\frac{\partial}
  {\partial y_{j}}
 $ becomes
 $
\frac{\delta}{\delta u^j}=\frac{\partial}{\partial u^j}
 $, therefore we have
$\partial_{\alpha}g_{11}= 0$.  Hence by a suitable choice of $t$, we
can assume
 \be\label{2.14}
 g_{11}= 1.
 \ee
 \epf
 \bp \label{spray}
 Let $(M,g)$ be a Finsler manifold admitting
 a  non-trivial
solution of (\ref{c-field}) . Then  the spray coefficients vanish.
 \ep
 \bpf We have $G^1= \frac{1}{2}\Gamma^{\ast}{}_{ij}^{1}\frac{du^i}{dt}\frac{du^j}{dt}$, therefore from (\ref{2.11})
  and (\ref{2.14}), we get $G^1=0$ and then from (\ref{G1}), we
  have
 \be \label{G}
 G^{i}=0.
 \ee
 \epf
As a consequence of the above proposition, $t$ may be regarded as
 the arc-length parameter of $t$-geodesics.
 \bt
Let $(M,g)$ be an $n$-dimensional  Finsler manifold admitting a
non-trivial solution of (\ref{c-field}), then $(M,g)$ is a
projectively flat Finsler manifold or is a direct product
$I\times\overline{M}$ of an open interval $I$ of the real line and
an $(n-1)$-dimensional Finsler manifold $\overline{M}$.
 \et
 \bpf
Let's consider the local coordinates
  $(t, u^2,...,u^n)$ on $M$, where $t$ is the arc-length parameter.
   Then the geodesic equation of $(M,g)$ becomes
   $\frac{d^2u^i}{dt^2}+\Gamma^{*{}i}_{jk}\frac{du^j}{dt}\frac{du^k}{dt}=0$ and from Proposition (\ref{spray}),
we get
 $\frac{d^2 u^{i}}{dt^2}=0$. If all  geodesics of $(M,g)$
are parameterized by  $t$ as arc-length then they are straight lines
and  by definition $(M,g)$ is a projectively flat Finsler manifold.
If not, a number of geodesics of $(M,g)$ should be parameterized by
$t$ and others by $(u^\alpha)$, then they will lie respectively on a
straight line and an hypersurface which is a $t$-level of $\rho$.
Therefore $(M,g)$ is a direct product $I\times\overline{M}$, where
$I$ is a real line and $\overline{M}$ is an $(n-1)$-dimensional
Finsler manifold, diffeomorphic to $t$-levels of $\rho$.
 \epf

 From (\ref{2.11}) and (\ref{2.14}), in  local  coordinates $(u^1=t, u^2,...,u^n)$,
the components of
  the Finsler metric tensor  $g$ is given by
 $$
(g_{ij})= \left(\begin{array}{lcl}
1\quad{0 \ldots\quad{ 0}}\\
0\quad{g_{22}\ldots g_{2n}}\\
\vdots \qquad{\ldots}\\
0\quad{g_{n2}\ldots g_{nn}}
\end{array}\right).
 $$
  \bl\label{2.19}
Let $(M,g)$ be a Finsler manifold and $\rho$  a non-trivial solution
of (\ref{c-field})  on $M$. Then the Finsler metric form of $M$ is
given by
  \be
ds^{2}=(du^{1})^{2}+ \rho'^{2}f_{\gamma\beta}du^{\gamma}du^{\beta},
  \ee
  where $f_{\gamma\beta}$ is a Finsler metric tensor
 on a $t$-level of $\rho$.
  \el
 \bpf
  Let $\overline{M}$ be a $t$-level of $\rho$. Then the unit vector field $\textbf{i}=\frac{grad \rho}{\|grad
  \rho\|}$, where $\|.\|=\sqrt{g(.,.)}$,
   is normal
to $\overline{M}$ at any point of $\overline{M}$ and the induced
metric tensor $\overline{g}_{\gamma\beta}$ of $\overline{M}$ is
given by $\overline{g}_{\gamma\beta}
=g_{ij}B_{\gamma}^{i}B_{\beta}^{j}$ where
$B_{\gamma}^{i}=\partial_{\gamma}u^{i}=\delta^i_\gamma$ (see [6]).
Therefore
 we have
 \be\label{g}
\overline{g}_{\alpha\beta}=g_{\alpha\beta}.
 \ee
 The $h$-second fundamental
form of $\overline{M}$ is defined by $
h_{\gamma\beta}:=(\nabla_{\gamma}^{H}B_{\beta}^{k})\textbf{i}_{k}$
and
 \be \label{2.7}
 \partial_{\beta}B_{\gamma}^{k}+
\Gamma^{\ast}_{ij}{}^{k}B_{\gamma}^{i}B_{\beta}^{j} -
\overline{\Gamma}_{\gamma\beta}^{\ast
\alpha}B_{\alpha}^{k}=h_{\gamma\beta}\textbf{i}^{k},
 \ee
 where $\overline{\Gamma}_{\gamma\beta}^{\ast \alpha}$ are Finsler
connection's coefficients  in
$(\overline{M},\overline{g}_{\gamma\beta})$. On the other hand,
$B_{\beta}^{k}\textbf{i}_{k}=0$, so we have
 \be \label{after h}
 h_{\gamma\beta}= -(\nabla _{j}^{H}
 \textbf{i}_{k})B_{\gamma}^{j}B_{\beta}^{k}.
 \ee
 Since the components of the unit vector field $\textbf{i}$ are
$\textbf{i}^{k}=\delta_{1}^{k}$, the equation (\ref{2.7}) for $k=1$
  reduces to
 \be \label{2.16}
 \Gamma^{\ast}{}^{1}_{\gamma\beta}=h_{\gamma\beta}.
  \ee
By h-covariant derivative of  $\rho_l= \textbf{i}_l{\|grad
  \rho\|}$ and using (\ref{c-field}), we have
 \be \label{*}
(\nabla_{k}^{H}{\|grad
  \rho\|})\textbf{i}_{l}+ {\|grad
  \rho\|}
(\nabla_{k}^{H}\textbf{i}_{l})=\phi g_{lk}.
 \ee
Contracting  with $\textbf{i}^{l}$, we have $ \nabla_{k}^{H}{\|grad
  \rho\|} = \phi \textbf{i}_{k}$, and by replacing  in (\ref{*}), we
   get
 \be \label{2.5}
 \nabla_{k}^{H}\textbf{i}_{l}=\frac{\phi}{\|grad \rho\|}(g_{lk}-
\textbf{i}_{k}\textbf{i}_{l}).
 \ee
Substituting (\ref{2.5}) into (\ref{after h}) , we have
 \be
 \label{2.8}
 h_{\gamma\beta}= h
g_{\gamma\beta},\quad{\quad{\textrm{where}\quad{
h=\frac{-\phi}{\|grad \rho\|}}}}.
 \ee
Using $\rho^l=g^{il}\rho_i$  and the fact that $\rho$ is a function
of $t$ alone, we have $\|\verb"grad" \rho\| =\rho'\neq0$, where
prime denotes the ordinary differentiation with respect to $t$. In
the same way we get
 \be
\label{2.17}\|\verb"grad" \rho\|=\rho', \quad{\phi=\rho'',
\quad{h=\frac{-\phi}{\|grad \rho\|}= \frac{-\rho''}{\rho'}}}.
 \ee
By using (\ref{2.8})  and replacing $h$  in (\ref{2.16}), we have
$\frac{\partial g_{\gamma\beta}}{\partial
t}=\frac{2\rho''}{\rho'}g_{\gamma\beta}$. Therefore  the components
$g_{\gamma\beta}$ are written in the form
 \be \label{2.18}
g_{\gamma\beta}=\rho'^{2}f_{\gamma\beta},
 \ee
where $f_{\gamma\beta}$ are functions of the $2(n-1)$ coordinates
$(u^{\alpha}, v^{\alpha})$. Since the metric tensor
$g_{\gamma\beta}$ is positive definite, so is the matrix
$(f_{\gamma\beta})$. Thus $(f_{\gamma\beta})$ can be regarded as
components of a Finsler metric tensor on $\overline{M}$.
 \epf
 \section{Curvature tensor of Cartan connection in adapted coordinates.}
Let $(M,g)$ be a Finsler manifold admitting a non-trivial solution
of Eq.(\ref{c-field}), we want to compute the components of Cartan
connection and its
 $h$-curvature tensor in terms of the adapted coordinates $(t, u^2,...,u^n)$.

 The coefficients $\rho'^{2}$ in (\ref{2.18}) are positive constants
in every $t$-levels of $\rho$. Therefore if $(N, f_{\gamma\beta})$
is an $(n-1)$-dimensional Finsler manifold diffeomorphic to a
$t$-level $(\overline{M}, \overline{g}_{\gamma\beta})$, then from
(\ref{g}), the Finsler manifold $(N, f_{\gamma\beta})$ and
$t$-levels
 neighboring $(\overline{M},
\overline{g}_{\gamma\beta})$ are locally homothetically
diffeomorphic to each other. Indeed, the connection coefficients
constructed from $f_{\gamma\beta}$ on $N$ have the same expression
as
 the connection
coefficients $\overline{\Gamma}_{\beta\gamma}^{\ast{}\alpha}$
constructed from the  induced metric $\overline{g}_{\gamma\beta}$ in
$\overline{M}$. Therefore $\Gamma_{ij}^{\ast {}k}$, the components
of  Cartan connection  on $(M,g)$, are given by
 \be \label{2.20}\begin{array}{ccc}
\Gamma^{\ast}{}_{11}^{1}=
\Gamma^{\ast}{}_{1\beta}^{1}=\Gamma^{\ast}{}_{\beta1}^{1}=\Gamma^{\ast}{}_{11}^{\alpha}=0,\\\\
\Gamma^{\ast}{}_{\gamma\beta}^{1}=-\frac{\rho''}{\rho'}g_{\gamma\beta}=-\rho'
\rho''f_{\gamma\beta},\\\\
\Gamma^{\ast}{}_{1\beta}^{\alpha}= \Gamma^{\ast}{}_{\beta
1}^{\alpha}=\frac{\rho''}{\rho'}\delta_{\beta}^{\alpha},\\\\
\Gamma^{\ast}{}_{\alpha\beta}^{\gamma}=\overline{\Gamma}^{\ast
\gamma}_{\alpha\beta},
 \end{array}
 \ee
where the last one comes from (\ref{2.7}), by replacing the index
$k$ with $\alpha$. The components of $h$-curvature tensor of Cartan
connection is given by $R^i_{\ hjk}=K^i_{\ hjk}+ C^{i}_{\
hr}R^{r}_{\ jk}$ where $K^i_{\ hjk}=\frac{\delta}{\delta
u^k}\Gamma_{hj}^{\ast}{}^{i}+\Gamma_{hj}^{\ast}{}^{r}\Gamma_{rk}^{\ast}{}^{i}-(j,k)$,
$R^{i}_{jk}=\frac{\delta}{\delta_{k}}G^{i}_{j}-(j,k)$ and
 $(j,k)$ denotes the interchange of indices $j$ and $k$. Here we have $R^{i}_{jk}=0$ and so $R^i_{\
hjk}=K^i_{\ hjk}$. Therefore $R^i_{\ hjk}$ of $(M,g_{ij})$, the
components of the $h$-curvature tensor,
 are given by
 \be \label{2.21}\begin{array}{ccc}
R^\alpha_{\ 1 \gamma 1}= - R^\alpha_{\ \gamma 1
1}=(\frac{\rho'''}{\rho'})\delta_{\gamma}^{\alpha},\\\\
R^1_{\ 1 \gamma \beta}= - R^ 1_{\ \gamma 1 \beta}=-\rho' \rho'''
f_{\gamma\beta},\\\\
R^\alpha_{\ \delta \gamma \beta}=\overline{R}^\alpha_{\ \delta
\gamma \beta}-(\rho'')^2(f_{\gamma\beta}\delta_{\delta}^{\alpha}-
f_{\delta\beta}\delta_{\gamma}^{\alpha}).
 \end{array}
 \ee

\section{ Critical points and their effects.}
Let $(M,g)$ be a Finsler manifold with Cartan connection and $\rho$
a solution of the differential equation (\ref{c-field}) on $M$. If
$g(\verb"grad" \rho,\verb"grad" \rho)=0$ in some points of $M$, then
$M$ possesses some interesting properties.

\bd
 A point $o$ of $(M,g)$ is called a critical point of $\rho$ if the
vector field $\verb"grad" \rho$ vanishes at  $o$, equivalently if $
\rho'(o)=0$. \ed

Let $\rho$ have a critical point $o$. We denote the distance
function from $o$ to an arbitrary point $p \in M$ by a
differentiable function $d(o,p):= inf L(\gamma)
 $, where $L(\gamma)=\int_a ^b
F(\gamma,\frac{d\gamma}{dt})dt$ and $\gamma:[a,b]\rightarrow M$ is a
piecewise $C^\infty$ curve with velocity
 $\frac{d\gamma}{dt}\in T_{\gamma(t)}M$ such that $\gamma(a)=o$ and
$\gamma(b)=p$.

We recall that, along any $t$-geodesic, Eq.(\ref{c-field}) reduces
to an ordinary differential equation
 \be\label{2.22}
\frac{d^{2}\rho}{dt^{2}}= \phi (\rho),
 \ee
where $\phi$ is a function of $\rho$ which is differentiable at
non-critical points.
 \bp
Let $(M,g)$ be a Finsler manifold  and $\rho$ a   solution of
Eq.(\ref{c-field})  on $M$ with a critical point $o$. If one of the
$t$-geodesics passes through $o$, then so do all of them.
 \ep
 \bpf
  Let $\overline{M}$ be a $t$-level of $\rho$, that is to say $\rho$ and $\rho'$  are constant on  $\overline{M}$ and
  there is no critical point on  $\overline{M}$. Let $p,q\in\overline{M}$ and denote by $\ell(p)$  and $\ell(q)$
   the $t$-geodesics
through  $p$ and $q$, respectively. The solution  of Eq.(\ref{2.22})
is given by the same function
   $\rho(t)$ on $\ell(p)$ and  $\ell(q)$, by uniqueness of solution of
ordinary differential equations. Hence if one of the
   $t$-geodesics
    passes through a critical
point at $t_0=0$, that is
  $\rho'(t_0)=0$, then so do all of them.
 \epf
 Moreover, from the uniqueness of solution of ordinary
differential equations (\ref{2.22}), if we denote a point at a
distance $t$ from $o$ on
   $\ell(p)$ by $p(t)$ and on $\ell(q)$ by $q(t)$, then
  the points  $p(t)$  and $q(t)$,
corresponding to the same value of $t$, are on the same
  $t$-level of $\rho$.
 \bp\label{isolate} Let $(M,g)$ be  a complete
Finsler manifold  and $\rho$  a solution of Eq. (\ref{c-field}) on
$M$ with a critical point $o$. Then $o$ is an isolated point and
$t$-levels of $\rho$
 are hyperspheres with center $o$.
 \ep
 \bpf
 Let  $U(o)$ be a geodesically connected neighborhood of $o$ and
$\overline{U}(o)$ be the closure of $U(o)$. If we denote by
$\overline{M}(p)$ a $t$-level of $\rho$ through a point $p$, then
$\overline{M}(p)$ is closed, so is the intersection $\overline{M}(p)
\cap \overline{U}(o)$, which contains its limit points. On the other
hand, there is no critical point on $\overline{M}(p)$, so $o$ is not
a limit point of $\overline{M}(p)$. Therefore the distance function
$d(o,p)$ between $o$ and the points of $\overline{M}(p)$ has a
minimum at an interior point of $U(o)$. Let $p_1$ be this interior
point and $t_1$ the minimum value. Thus the $t$-geodesic which joins
$o$ to $p$ is  $\ell(p_1)$. We denote the geodesic hypersphere with
center $o$ and radius $t_1$ by $S^{n-1}(o,t_1)$. Since
$\overline{M}(p) \cap S^{n-1}(o,t_1)$ is both closed and open  in
$\overline{M}(p)$, by taking into account the connectedness of
$\overline{M}(p)$, it becomes evident that the geodesic hypersphere
$S^{n-1}(o,t_1)$ coincides with the $t$-level $\overline{M}(p)$ and
the proposition is proved.
 \epf
 As a consequence of the above proposition, one can show easily  that the
number of critical points of $\rho$ is not more than two.
\section{Spaces of constant curvature.}
 Let $P(v,X)\subset T_{u}(M)$ be a 2-plane generated by $v$
and $X\in T_{u}(M) $, where $(u,v)$ is the line element of $TM$. The
\emph{sectional curvature} with respect to $P$ is given by
$$K(u,v,X):=\frac{g(R(X,v)v,X)}{g(X,X)g(v,v)-g(X,v)^{2}}.$$
If $K$ is independent of $X$, then $(M,g)$ is called \emph{space of}
\emph{scalar curvature}. If $K$ has no dependence on $u$ or $v$,
then the Finsler manifold is said to have \emph{constant curvature},
see [1] or [3].
 \bl \label{prop}
  Let $(M,g)$ be  an n-dimensional Finsler manifold which admits a
solution
 $\rho$  of Eq. (\ref{c-field}) with one critical
point. Then the $(n-1)$-dimensional  manifold $\overline{M}$ with
Finsler metric form $\overline{ds}^{2}= f_{\gamma\beta}du^{\gamma}
du^{\beta}$ as defined in (\ref{2.18}) is a space of positive
constant curvature.
 \el
\begin{proof} Let $o$ be a critical
point of $\rho$. Proposition \ref{isolate} implies that, there is a
geodesically connected neighborhood $U(o)$ of $o$, for which only
the point  $o$ is critical. Hence from Proposition \ref{2.19}, the
Finsler metric form in $U(o)$ becomes
 $$ds^{2}=(du^{1})^{2}+
\rho'^{2}\overline{ds}^{2},$$
 where $\overline{ds}^{2}=f_{\gamma\beta}du^{\gamma} du^{\beta}$ is the
Finsler metric form of an $(n-1)$-dimensional manifold
$\overline{M}$ diffeomorphic to the $t$-levels of $\rho$. Now we can
consider $\overline{R}^\alpha_{\ \delta\gamma\beta}$ in
(\ref{2.21}), as  the components of the $h$-curvature tensor of
$\overline{M}$. Therefore the norm of the $h$-curvature tensor
$R^h_{\ ijk}$ with respect to the metric tensor $g$ is given by
 $$
\|R^h_{\ ijk}\|_{g}^{2}=R^h_{\ ijk}R_h^{\
ijk}=\frac{1}{\rho'^4}\|\overline{R}^\alpha_{\ \delta\gamma\beta}-
\rho''^{2}(f_{\gamma\beta}\delta_{\delta}^{\alpha}-
f_{\delta\beta}\delta_{\gamma}^{\alpha})\|^{2}_{f} +
4n(\frac{\rho'''}{\rho'})^{2}.
 $$
By definition, in a critical point $o$ at $t=0$ we have $\rho'(0)=
0$. From $\rho''(0)= \phi(o)$ in (\ref{2.17}) and the fact that
$f_{\gamma\beta}$, the components of the metric tensor $f$ and
$h$-curvature tensor $\overline{R}^\alpha_{\delta\gamma\beta}$ of
$\overline{M}$ are independent of $t$, we can conclude that the
above equation, as $t$ tends to zero, becomes
 $$
\overline{R}^\alpha_{\delta\gamma\beta}=\phi(o)^{2}(f_{\gamma\beta}\delta_{\delta}^{\alpha}-
f_{\delta\beta}\delta_{\gamma}^{\alpha}).
 $$
By means of Proposition \ref{isolate}, the $t$-levels of $\rho$ are
hyperspheres with center $o$ and therefore $\phi(o)\neq 0$. Hence
$\overline{M}$ has positive constant sectional curvature
$\overline{K}=\phi(o)^{2}$.
\end{proof}
As it is mentioned on the proof  of the above lemma, $\phi(o)$ does
not vanish in this case and we have
 \be \label{2.23}
\lim_{t\rightarrow0}\frac{\rho'(t)}{t}=\rho''(0)=\phi(o)\neq 0.
 \ee
Therefore $\rho'(t)$ and $t$ are of the same order.

In  case there is a solution of (\ref{c-field}) with two critical
points, and the Finsler manifold $(M,g)$ is compact, an extension of
Milnor Theorem [16] implies that  $M$ is homeomorphic to an
n-sphere. Thus we have the following proposition.
 \bp
Let $(M,g)$ (dim $M>2$) be a simply connected and compact Finsler
manifold which admits a solution $\rho$ of (\ref{c-field}) with two
critical points, then $M$ is homeomorphic to an n-sphere.
 \ep

\section{A classification of complete Finsler manifolds.}

 Here we summarize the  above results on the existence of solutions of Eq.(\ref{c-field}).
 \bt\label{main}
Let $(M,g)$ be a connected complete Finsler manifold of dimension
$n\geq 2$. If $M$ admits a non-trivial solution of $
\nabla^{H}\nabla^{H}\rho= \phi g$, where $\nabla^{H}$ is the Cartan
horizontal covariant derivative,
 then  depending on  the number of critical points of $\rho$, i.e. zero, one or two respectively,
  it is conformal to\\
 \textbf{(a)} A direct product $J\times\overline{M}$ of an open
interval $J$ of the real line and an $(n-1)$-dimensional complete
Finsler manifold $\overline{M}$.\\ \textbf{(b)} An n-dimensional Euclidean space.\\
 \textbf{(c)} An n-dimensional  unit sphere in an Euclidean
 space.
 \et
 \begin{proof}Let $\overline{M_{0}}$ be a $t$-level of $\rho$ and $\overline{M}$ an $(n-1)$-dimensional
Finsler manifold having the metric form
$\overline{ds}^{2}=f_{\gamma\beta} du^{\gamma} du^{\beta}$ defined
as in (\ref{2.18}), which coincides with $\overline{M_{0}}$ as a set
of points. First of all we note that in a complete manifold $M$,
$\overline{M_{0}}$ is complete with respect to the induced Finsler
metric. In fact, the distance between points of $\overline{M_{0}}$
with respect to the induced metric is not shorter than that of $M$,
and hence a Cauchy sequence of points in $\overline{M_{0}}$ is also
a Cauchy sequence in $M$. Since $M$ is complete and
$\overline{M_{0}}$ is closed in $M$, every Cauchy sequence has its
limiting point in $\overline{M_{0}}$, hence $\overline{M_{0}}$ is
complete. By means of (\ref{2.18}), the induced  tensor metric  of
$\overline{M_{0}}$ is proportional to the tensor metric of
$\overline{M}$, so they are homothetic to each other. Therefore
$\overline{M}$ is also complete.

\textbf{case (a)} If one of the $t$-geodesics orthogonal to
$\overline{M_{0}}$ has no critical point, then non of them has.
Since $M$ is complete, $t$-geodesics are extendable to whole
interval $I=(-\infty , +\infty)$ of the arc-length $t$. So we can
define the following map
$$\nu: I \times \overline{M}\rightarrow M,$$
$$\quad{(t,p)\rightarrow p(t)},$$
where the point $p(t)$ corresponding to the value $t \in I$ lies on
a  $t$-level of $\rho$ and we have $\nu(0, \overline{M})=
\overline{M}_{0}$. On the other hand, $\overline{M}$ is
diffeomorphic  to the $t$-levels, therefore the map $\nu$ is a
diffeomorphism of $I\times\overline{M}$ into $M$. Since $M$ is
connected, any point $q$ of $M$ is joined to a point of
$\overline{M}_{0}$ by a curve $\mathcal{C}$. By extending local
diffeomorphism among $t$-levels of $\rho$ through points of
$\mathcal{C}$, we can see that $q$ is an image of a point $(t,p)$ of
$I \times \overline{M}$. Thus the map $\nu$ is a diffeomorphism of
$I\times\overline{M}$ onto $M$. Therefore from the proof of Lemma
\ref{prop}, the metric form of $M$ is expressed as
 \be \label{2.24}
 ds^{2}= dt^{2}+
(\rho')^{2}\overline{ds}^{2},
 \ee
where $\overline{ds}^{2}$ is the metric form of $\overline{M}$.

After a reparametrization of $t$-geodesics such that $\rho'>0$, we
define a parameter $r$ by
 \be
\label{2.25} r(t)=\int_{0}^{t}\frac{1}{\rho'}dt\quad{\quad{t\in I}}.
 \ee
$r(t)$ is an increasing monotone function of $t$. Let
 $$
r_{1} = \lim_{t\rightarrow-\infty}r(t),    \quad{\textrm{and}}
\quad{r_{2}= \lim_{t\rightarrow+\infty}r(t)},
 $$
where $r_{1}$ and $ r_{2}$ may be infinite, and let $J$ be the
interval $(r_{1}, r_{2})$. Now because $dr=\frac{1}{\rho'}dt$, we
can write the metric form (\ref{2.24}) as follows
 \be
 \label{2.26}
ds^{2}=(\rho')^{2} (dr^{2} + \overline{ds}^{2} ) , \quad{\quad{r\in
J}}.
 \ee
Thus $M$ is conformal to the direct product $J\times \overline{M}$.

\textbf{case (b)}
 If one of the $t$-geodesics issuing from points of
$\overline{M_{0}}$ has a critical point $o$ at a distance $t_{0}$
from $\overline{M_{0}}$ and no critical point in the opposite
direction, then  all such curves have the same behavior, and $o$ is
the only critical point of $\rho$. Let's parameterize $t$-geodesics
by arc-length $t$ measured from $o$ and put $I=(0,\infty)$. The map
$\nu: I \times \overline{M}\rightarrow M$ defined  as in the case
(a), is a diffeomorphism of $I \times \overline{M}$ onto the open
set $M \setminus \{o\}$. By Lemma \ref{prop}, $\overline{M}$ with
the metric form $\overline{ds}^{2}$ is a space of positive constant
curvature $(\phi(o))^{2}$. If we suppose that $\phi(o)>0$ in $I$ and
put $\overline{c}= \phi(o)$, then we can define a parameter $r$ for
the $t$-geodesics by
 \be \label{2.27}
r(t) = e^{\overline{c}\int_{t_{0}}^{t} \frac{dt}{\rho'}},\quad
{\quad{ t\in I}},
 \ee
which is an increasing monotone function of $t$.
 From (\ref{2.23}),
we see that $\rho'(t)$ is of the same
 order as $t$, when $t$ tends to zero, hence $$\lim_{t\rightarrow
 0}r(t)= e^{\overline{c}(\lim_{t\rightarrow
 0}\int\frac{dt}{t})}=e^{\overline{c}(-\infty)}
 =0.$$
 If we put $r_{2}= \lim_{t\rightarrow \infty}r(t) \leq
 \infty$, then the parameter $r$ varies in the interval $[0,r_{2}
 )$ as $t$ varies in $[0,+\infty)$. Since
 $$
\frac{dr}{r}=\frac{\bar{c}dt}{\rho'(t)},
 $$
the metric form of $M$ is equal to
 $$
ds^{2}= (\frac{\rho'(t)}{r(t)\overline{c}})^{2}[dr^{2}+
r^{2}\overline{c}^{2}\overline{ds}^{2}],\quad{\quad{0<r<r_{2}}},
 $$
in $M \setminus \{o\}$. The expression in the brackets is the polar
form of an $n$-dimensional Euclidean metric. On the other hand,
taking into account  orders of $\rho'(t)$ and $r(t)$, we can see
that the coefficient $\frac{\rho'(t)}{r(t)}$ is not equal to zero
but it is differentiable at $o$. Thus $M$ is conformal to the
Euclidean ball of radius $r$. Since $r_2$ can be increased without
bound, $M$ is conformal to a flat space and by definition, it is
conformally flat.

\textbf{case (c)}
 If one of the $t$-geodesics issuing from points of
$\overline{M_{0}}$ has two critical points $o$ and $o'$ in opposite
directions at distances $t_{1}$ and $t_{2}$ respectively, then so do
all such curves, and only the points $o$ and $o'$ are critical in
$M$. We parameterize the $t$-geodesics by  the arc-length $t$
measured from $o$. Let the distance from $o$ to $o'$ be equal to
$2t_{0}$, and put $I=(0 , 2t_{0})$. Moreover let $\overline{M_{0}}$
be the $t$-level of $\rho$ corresponding to $t_{0}$.
 The map
$\upsilon : I\times \overline{M} \longrightarrow M$, defined  as in
the case $(a)$, is a diffeomorphism of $I\times \overline{M}$ onto
the open set $M\setminus{\{o,o'\}}$. The metric form of $M$ is
written as in (\ref{2.24}), where from the Lemma \ref{prop}, the
manifold $\overline{M}$ with metric form $\overline{ds}^{2}$ is a
space of positive constant curvature
$\overline{k}=\phi(o)^{2}=\phi(o')^{2}$ . Let $\rho'(t)>0$ on the
interval $I$, that is  $\rho$ is an increasing function on $I$, so
we have $\rho''(0)>0$ and $\rho''(2t_{0})<0$. Now we can define a
parameter $ \theta$ by
 \begin{equation*}\label{2.28}
\theta(t)= 2 \arctan\exp\overline{c} \int_{t_{0}}^{t}
\frac{dt}{\rho'(t)},\quad {\quad{ t\in I}},
 \end{equation*}
 where $\overline{c}= \phi(o)=-\phi(o')$. $\theta$
is an increasing monotone  function of $t$ and we have
 $$
\lim_{t\longrightarrow 0}\theta(t)=
0,\quad{\theta(t_{0})=\frac{\pi}{2}},\quad{ \lim_{t\longrightarrow
2t_{0}}\theta(t)= \pi}.
 $$
Hence $\theta$ varies in the closed interval $[0,\pi]$ as $t$ varies
in $[0,2t_{0}]$. Since we have
 $$
\frac{d\theta}{\sin\theta}=\frac{\bar{c}dt}{\rho'(t)},
 $$
the metric form of $M$ is equal to
 \be \label{2.29}
ds^{2}=(\frac{\rho'(t)}{\overline{c}\sin\theta(t)})^{2}[d\theta^{2}+(\sin\theta)^{2}
\overline{c}^{2}\overline{ds}^{2}],\quad{\quad{0<\theta<\pi}},
 \ee
in $M\setminus{\{o,o'\}}$. The expression in the brackets is the
  polar form of the Finsler metric on an n-sphere [18].
Taking into account orders of $\rho'(u)$ and $\sin\theta(u)$, we can
verify that the factor $\rho'(u)/\sin\theta(u)$ is not equal to zero
but it is differentiable at both the critical points $o$ and $o'$.
Therefore $(M,g)$ is conformal to a Finsler metric on an n-sphere.
\end{proof}
\section{ Example of  Finsler metrics with
positive constant curvature.} Describing the Finsler metrics of
constant flag curvature is one of the fundamental problems in
Finsler geometry.  Historically, the first set of non-Riemannian
Finsler metrics of constant flag curvature are the Hilbert-Klein
metric and the Funk metric on a strongly convex domain.  In 1963,
Funk [12] completely determined the local structure of
two-dimensional projectively flat Finsler metrics with constant flag
curvature. The Funk metric is positively complete and non-reversible
with $K= -\frac{1}{4}$ and the Hilbert-Klein metric is complete and
reversible with $K= -1$. Both of them  are locally projectively
flat.
 Yasuda and Shimada [23] in 1977, classified Randers
metrics of constant flag curvature, which has been rectified and
completed in a joint work of D. Bao, C. Robles and Z. Shen [4] in
2004. Akbar-Zadeh [2] in 1988, proved that a closed Finsler manifold
with constant flag curvature $K$ is locally Minkowskian if $K= 0$,
and Riemannian if $K= -1$. In the case $K= 1$, Shen [19] asserts
that the Finsler manifold must be diffeomorphic to sphere, provided
that it is simply connected. Then, Bryant [7] showed that up to
diffeomorphism, there is exactly a 2-parameter family of locally
projectively flat Finsler metrics on $S^2$ with K= 1 and the only
reversible one is the standard Riemannian metric. He has also
extended his construction to higher dimensional spheres . These
Bryant's examples are projectively flat and none of them is of
Randers type. In 2000, Bao and Shen [5] constructed a family of
non-projectively flat Finsler metrics on $S^3$ with $K=1$, using the
Lie group structure of $S^3$. They also produced, for each constant
$K>1$, an explicit example of a compact boundaryless
non-projectively flat Randers space with constant positive flag
curvature $K$.

 Here, based on a Finsler
metric of positive constant curvature on certain hypersurfaces of
the Finsler space  $(M,g)$, we find some conditions for $g$ to be a
Finsler metric of positive constant curvature. This constructed
Finsler metric is not necessarily of Randers type nor projectively
flat.
 Without loss  of generality,  we have so far considered
some hypersurfaces of $I\!\!R^{n+1}$.
 \bp
Let $\overline{M}$ be a regular hypersurface of $I\!\!R^{n+1}$
defined by $\rho=$constant. The scalar function $\rho$ is a solution
of
 $\nabla^{H}\nabla^{H}\rho=K^2 \rho g$  on the Finsler space $(I\!\!R^{n+1}, g)$ and $K$ is a constant number.
  If there exists  a Finsler
metric $\overline{g}$ defined on $\overline{M}$ with positive
constant curvature, then  $g$ is a Finsler metric of this kind on
$I\!\!R^{n+1}$.
 \ep
 \bpf
 Let  $\rho$  be a differentiable scalar function on the Finsler space  $(I\!\!R^{n+1}, g_{ij})$
 satisfying
 \be\label{K}
 \nabla^{H}_i\nabla^{H}_j\rho=K^2 \rho g_{ij},
 \ee
where $\nabla^{H}$ denotes the Cartan $h$-covariant derivative, $K$
is a positive constant, $(u^i, v^i)$  is the  local coordinate
system on $TI\!\!R^{n+1}$  and $i,j$ run over the range $1,...,n+1$.
Along the geodesics of $(I\!\!R^{n+1},g)$  with arc-length parameter
$t=u^{n+1}$, equation (\ref{K}) reduces to $\frac{d^2
\rho}{(dt)^2}=K^2 \rho$. The general solution of this differential
equation is given by $\rho (t)=a \cos (K t)+b$, where $a$ and $b$
are constants. Differentiating with respect to $t$ gives $\rho '(t)=
-aK \sin K t$. Therefore, there are two critical points
corresponding to $t=0$ and $t=\frac{\pi}{K}$. Let $\overline{M}$ be
the  hypersurface of $I\!\!R^{n+1}$  defined by $\rho=$constant and
$\overline{g}$ a Finsler metric with positive constant curvature
$\overline{K}$ on $\overline{M}$. If we put $a=\frac{-1}{K}$, then
from (\ref{2.29}) the Finsler structure of $I\!\!R^{n+1}$ becomes
 \be\label{meter}
ds^2:=g_{ij}du^idu^j=(dt)^2+(\sin K t)^2\overline{K}^2
\overline{g}_{\alpha\beta}du^\alpha du^\beta,
 \ee
 where  $\alpha,\beta =1,..., n$. This is the  polar form of a
 Finsler
 metric on an $n$-sphere in $I\!\!R^{n+1}$ with positive  constant curvature $K$ [18].
 \epf
The same method can be applied to construct a Finsler metric with
positive constant curvature on an $(n+1)-$dimensional complete
Finsler manifold $(M,g)$, starting from a Finsler metric of this
kind on a hypersurface of $M$.

Behroz BIDABAD \\
{\small E-mail address: bidabad@aut.ac.ir}\\
Azam ASANJARANI \\
{\small E-mail address: asanjarani@aut.ac.ir}

\end{document}